\documentclass[12pt,a4 paper]{article}
\usepackage{amssymb}
\usepackage{amsmath}
\usepackage{amsthm}
\usepackage{amsopn}
\def \bP{{\mathbf P}}

\let\goth\mathfrak   

\def \zs#1{_{\lower 3pt \hbox{$\scriptstyle#1$}}}

\def \Ex{\mathbf E}
\theoremstyle{plain}
\newtheorem{theorem}{Theorem}

\newtheorem{corollary}{Corollary}

\newcommand{\abs}[1]{\left\vert#1\right\vert}

\def \bE{{\mathbf E}}

\begin{document}
\title{Power Loss for Inhomogeneous Poisson Processes}
\author{Fazli Kh.\\
{\small  Department of Mathematics, }\\
{\small University of Kurdistan, Sanandaj}\\
{\small khfazli@uok.ac.ir}}
 \maketitle

\begin{abstract}
In this work, based on a realization of an inhomogeneous Poisson
process whose intensity function depends on a real unknown
parameter, we consider a simple hypothesis against a sequence of
close (contiguous) alternatives. Under certain regularity conditions
we obtain the power loss of the score test with respect to the
Neyman-Pearson test. The power loss measures the performance of a
second order efficient test by the help of third order asymptotic
properties of the problem under consideration.
\end{abstract}



\section{Introduction}
Let $X^{(n)}$ be a realization of a nonhomogeneous Poisson process
observed on some increasing subsets $\mathbb{A}_n,\; n=1,2,\ldots$
of $d$ dimensional Euclidian space $\mathbb{R}^d$ with intensity
function $S\left(\vartheta,x\right), x\in \mathbb{A}_{n}$ depending
on one-dimensional parameter $\vartheta \in \Theta$. Based on
$X^{(n)}$ we want to test the hypotheses
\begin{align*}
&{\cal H}_{0}:\; \vartheta=\vartheta_{0}\\
&{\cal H}_{1}:\; \vartheta>\vartheta_{0},
\end{align*}
where $\vartheta_0$ is a given value in the parameter space $\Theta$
. Let us fix some $\alpha\in(0,1)$ and define the class ${\cal
K}^{(n)}_{\alpha}$ of tests at level $1-\alpha$ (size $\alpha$),
i.e.,
$$
{\cal K}^{(n)}_{\alpha}=\left\{\phi_n:\quad
\Ex_{\vartheta_{0}}\phi_n\left(X^{\left(n\right)}\right)=\alpha
\right\},
$$
where $\Ex_{\vartheta}$  denotes the mathematical expectation with
respect to the probability measure $\bP^{(n)}_{\vartheta}$. With
fixed $n$, generally speaking, there is no uniformly most powerful
test in  ${\cal K}^{(n)}_{\alpha}$ (see \cite{Pfan2}). Therefore we
turn to the asymptotic approach and introduce the class ${\cal
K}'_{\alpha}$ of sequence of tests of asymptotic level $1-\alpha$,
i.e.,
$$
{\cal K}'_{\alpha} =\left\{\{\phi_n\}:\quad \lim_{n\rightarrow
\infty }\Ex_{\vartheta_{0}}
\phi_n\left(X^{\left(n\right)}\right)=\alpha \right\}.
$$
It is well known that if $n\rightarrow\infty$ for any given value
$\vartheta$ of the alternative the power of any reasonable
(consistent) test tends to 1 ( see \cite{Pfan1}, \cite{Ben}). In
order to compare the different tests we use the Pitman's approach
(see \cite{Pitm}, \cite{Noe}) where instead of a fixed alternative
$\vartheta$ we consider a sequence of so-called
 {\sl local alternatives }({\sl close} or {\sl contiguous
alternatives}) which  converges to $\vartheta_0$ with a certain rate
and hence it is difficult to distinguish between the null hypothesis
and alternative. More precisely, let $\{\varphi_n\}$ be a sequence
of nonnegative numbers which converges to zero with such a rate that
the likelihood ratio
$$
Z_n(u)=\frac{{\rm d}\bP^{(n)}_{\vartheta_0+\varphi_{n}u}}{{\rm
d}\bP^{(n)}_{\vartheta_{0}}}\left(X^{\left(n\right)}\right)
$$
has a nondegenerate limit for any $u$ with $\vartheta_0+\varphi_n
u\in \Theta$. By the Neyman-Pearson lemma the most powerful test for
$\mathcal{H}_0:\;\vartheta=\vartheta_0$ against the local
alternative $\mathcal{H}_u:\;\vartheta=\vartheta_0+\varphi_n u$ with
$u>0$, is given by
\begin{equation*}
\tilde{\phi}_n\left(X^{\left(n\right)}\right)=\begin{cases}1,& \,
{\rm if}\,\, \Lambda_n(u)>b_n(u)\\
0,& \, {\rm if}\,\, \Lambda_n(u)<b_n(u)
\end{cases}
\end{equation*}
where $\Lambda_n(u)=\ln Z_n(u)$ and the constant $b_n(u)$ together
with the contribution of the randomized part provide the size
$\bE_{\vartheta_0}\tilde{\phi}_n\left(X^{\left(n\right)}\right)=\alpha$.
The power of $\tilde{\phi}_n$ as a function of $u$ is called the
{\sl envelope power function }. For any fixed $n$ it is the supremum
of the power at the local alternative
$\vartheta_u=\vartheta_0+\varphi_n u$ over all the tests at level
$1-\alpha$, i.e.
$$
\bE_{\vartheta_u}\tilde{\phi}_n\left(X^{\left(n\right)}\right)=\sup_{\phi_n\in
\mathcal{K}_\alpha^{(n)}}
\bE_{\vartheta_u}\phi_n\left(X^{\left(n\right)}\right).
$$
Notice that $\tilde{\phi}_n$ is not a test for the main hypotheses $
\mathcal{H}_0$ and $ \mathcal{H}_1$ because it depends on the
parameter $u$.

Let us introduce the score statistic
\begin{align*}
\Delta_{n}\left(\vartheta_{0}\right)=\varphi_{n}
\int_{\mathbb{A}_{n}}\frac{\dot{S}\left(\vartheta_{0}, x
\right)}{S\left(\vartheta_{0}, x \right)}\,\pi^{(n)}\left({\rm
d}x\right),\quad \varphi_n^{-2}=\int_{\mathbb{A}_{n}}
\frac{\dot{S}\left(\vartheta_0, x \right)^2}{S\left(\vartheta_0, x
\right)}\;{\rm d}x
\end{align*}
where the normalizing factor $\varphi_n=I_n(\vartheta_0)^{-1/2}$ is
the inverse square root of the {\sl Fisher information}
$$
I_n(\vartheta_0)=\int_{\mathbb{A}_{n}}
\frac{\dot{S}\left(\vartheta_0, x \right)^2}{S\left(\vartheta_0, x
\right)}\;{\rm d}x
$$
at the point $\vartheta_0$. Here $\pi^{(n)}\left({\rm
d}x\right)=X^{(n)}\left({\rm d}x\right)-S\left(\vartheta_0, x
\right){\rm d}x$\, is the centered Poisson process (for the
definition of a stochastic integral w.r.t. a Poisson process see the
next section) and $\dot{S}\left(\vartheta, x \right)$ denotes the
derivative of $S\left(\vartheta, x \right)$ with respect to
$\vartheta$. Based on $\Delta_{n}\left(\vartheta_{0}\right)$ we
introduce the score test
\begin{equation*}
\bar{\phi}_n\left(X^{\left(n\right)}\right)=\begin{cases}1,& \,
{\rm if}\,\, \Delta_{n}\left(\vartheta_{0}\right)>z_\alpha\\
0,& \, {\rm if}\,\, \Delta_{n}\left(\vartheta_{0}\right)\leq
z_\alpha
\end{cases}
\end{equation*}
where $z_\alpha$ is $1-\alpha$ quantile of standard Gaussian law,
{\sl i.e.}, $\bP\left\{\zeta > z_\alpha\right\}=\alpha$ and
$\zeta\sim\mathcal{N}(0,1)$. It is well known that if the family
$\left\{\bP^{(n)}_{\vartheta},\, \vartheta\in \Theta \right\}$ of
distributions is {\sl locally asymptotically normal} (LAN) at the
point $\vartheta_{0}$, then the test $\bar{\phi}_n\in
\mathcal{K}'_\alpha$ is {\sl locally asymptotically uniformly most
powerful } (LAUMP) (or {\sl first order efficient}), {\sl i.e.}, for
any $K>0$
$$
\sup_{0\leq u \leq
K}\abs{\bE_{\vartheta_u}\tilde{\phi}_n\left(X^{\left(n\right)}\right)-
\bE_{\vartheta_u}\bar{\phi}_n\left(X^{\left(n\right)}\right)}=o(1)
$$
as $n\rightarrow\infty$ (see \cite{Rou}). Moreover the power
function of $\bar{\phi}_n$ at $\vartheta_u$ admits the
representation
\begin{equation}
\label{power1}
\bE_{\vartheta_u}\bar{\phi}_n\left(X^{\left(n\right)}\right)=\bP\left\{\zeta
> z_{\alpha }- u \right\}+ o\left(1\right)
\end{equation}
for any $u> 0$, where $\zeta\sim\mathcal{N}(0,1)$. For $n$ large,
hence the power of $\bar{\phi}_n$ approximates the envelope power
function up to order $o(1)$. A refinement of \eqref{power1} is given
in \eqref{power21}. For a family $\left\{\bP^{(n)}_{\vartheta},\,
\vartheta\in \Theta \right\}$ of distributions related to a Poisson
process with intensity functions
$\left\{S\left(\vartheta,\cdot\right),\, \vartheta\in \Theta
\right\}$ the conditions of LAN for multidimensional parameter
$\vartheta$ are obtained by Yu. A. Kutoyants, \cite{Kut98}. The
first order efficiency of $\bar{\phi}_n$ follows from the LAN
representation which implies in turn the asymptotic normalities of
$\Delta_n(\vartheta_0)$ and $\Lambda_n(u)$ under both
$\mathcal{H}_0$ and $\mathcal{H}_u$. Therefore the refinement of the
central limit theorem, by taking into  account one term after the
Gaussian term, improves the situation. This can be done by the help
of Edgeworth type expansion of the distribution function of the
stochastic integral $\Delta_n(\vartheta_0)$ and $\Lambda_n(u)$ under
$\mathcal{H}_0$ and $\mathcal{H}_u$. Under certain regularity
conditions related to second order asymptotic properties of the
family $\left\{\bP^{(n)}_{\vartheta},\, \vartheta\in \Theta
\right\}$, we can construct a {\sl second order efficient} test,
{\sl i.e.,} a test ${\phi}_n^{*}$ such that for any $K>0$ 
\begin{equation}
\label{sec}\sup_{0\leq u\leq
K}\abs{\bE_{\vartheta_u}\tilde{\phi}_n\left(X^{\left(n\right)}\right)-
\bE_{\vartheta_u}\phi_n^*\left(X^{\left(n\right)}\right)}=O(\varepsilon_n^2),
\end{equation}
for some sequence  $\varepsilon_n\rightarrow 0$ (see \cite{Fa1}).
Furthermore the probability of the first type error of
${\phi}_n^{*}$ is given by
$\bE_{\vartheta_0}\phi^{*}_{n}\left(X^{(n)}\right)=
\alpha+O(\varepsilon_n^2)$. From \eqref{sec} it follows that the
power function of $\phi^{*}_{n}$ (a second order efficient test,
generally) approximates the envelope power function up to order
$O(\varepsilon_n^2)$ and hence it works as good as $\tilde{\phi}_n$
up to this order. Second order efficiency of $\phi^{*}_{n}$ is
related to the fact that the two first terms in the Edgeworth
expansions of the distributions functions of $\Delta_n(\vartheta_0)$
and $\Lambda_n(u)$ under the local alternative are equal up to the
order $O(\varepsilon_n^2)$. Hence to measure the performance of a
second order efficient test, and especially $\phi^{*}_{n}$, it is
natural to consider the {\sl power loss} of $\phi^{*}_{n}$ with
respect to the most powerful test
$\tilde{\phi}_{n}$, which is defined by 
\begin{equation}
r(u)=\lim_{n\rightarrow\infty}
\varepsilon_{n}^{-2}\left(\bE_{\vartheta_u}\tilde{\phi}_n\left(X^{\left(n\right)}\right)-
\bE_{\vartheta_u}\phi_n^*\left(X^{\left(n\right)}\right)\right),
\end{equation}
for $u>0$. This requires to take into account higher order terms in
the Edgeworth expansions of the distribution functions. See
\cite{Ben}, chapter 3, for a general theorem and the power loss
results for the tests based on $L-, R-$ and $U-$ statistics in the
i.i.d. case. The main object of this work is to obtain the power
loss of the score test $\phi_n^*$ based on a realization
$X^{\left(n\right)}$ of a nonhomogeneous Poisson process with
intensity function $S\left(\vartheta,x\right), x\in \mathbb{A}_{n}$
and to give the explicit representation of $r(u)$.

\section{Preliminaries}
Let us remind several facts from a Poisson process. A Poisson
process $X^{(n)}$ is a random point measure which on the set
$\mathbb{B}\subset \mathbb{A}_{n}$ can be written as
$$
X^{(n)}\left(\mathbb{B}\right)=\sum_{x_i\in \mathbb{A}_{n}}^{}\;
\chi\zs{\left\{x_i\in \mathbb{B}\right\}} ,
$$
where $\left\{x_i\right\}$ are the {\sl events} (random points) of
the Poisson process and $\chi\zs{\left\{\mathbb{D}\right\}} $ is the
indicator function of the event $\mathbb{D}$. The Poisson process
with (parametric) intensity function
$S(\vartheta,x),\;x\in\mathbb{A}_n$ (with respect to Lebesgue
measure) is entirely defined by the following two conditions:
\begin{itemize}
\item for any collection of disjoint sets
$\mathbb{B}_{1},\ldots,\mathbb{B}_{m}\subseteq\mathbb{A}_{n}$ the
random variables
$X^{(n)}(\mathbb{B}_{1}),\ldots,X^{(n)}(\mathbb{B}_{m})$ are
independent,
\item  the random variable $X^{(n)}(\mathbb{B})$ for any
$\mathbb{B}\subseteq\mathbb{A}_{n}$ has Poisson distribution  with
parameter
$\Lambda_{\vartheta}^{(n)}(\mathbb{B})=\int_{\mathbb{B}}S(\vartheta,x)\;{\rm
d}x$.
\end{itemize}

This and other definitions of the spatial Poisson processes as well
as their properties and examples can be found in  many books devoted
to  point processes (see, e.g., Daley and Vere-Jones \cite{DaVer},
Krickeberg \cite{Krik2}, Reiss \cite{Reiss}, Ripley \cite{Rip},
Snyder and Miller \cite{SM}). The spatial Poisson processes are
widely used in many fields. In particular this is a good
mathematical model for the problems of image restoration when the
optical signal is weak and statistics of photons is well described
by an inhomogeneous Poisson process \cite{SM}.

By the definition, $X^{(n)}$ is a random element of the set
$\mathcal{M}_{0}^{(n)}$, containing all integer valued measures
defined on the set $\mathbb{A}_{n}.$ Let
$\goth{B}(\mathcal{M}_{0}^{(n)})$ be the smallest $\sigma-$field
with respect to which all the mappings:
$$\Pi_{\mathbb{B}}:\mathcal{M}_{0}^{(n)}\rightarrow\{0,1,2,\ldots,\infty\},
\hspace*{1cm}\Pi_{\mathbb{B}}(X^{(n)})=X^{(n)}(\mathbb{B}),
\hspace*{0.5cm}\mathbb{B}\in\mathcal{A}_{n},
$$
are measurable. Here  $\mathcal{A}_{n}$ is the $\sigma-$field of
Borel subsets of  $\mathbb{A}_n$.  Let $\bP^{(n)}_{\vartheta}$
denote the probability law induced by the random element
(realization) $X^{(n)}$ of a Poisson process with intensity function
$S(\vartheta,x),\;x\in\mathbb{A}_n$ on the measurable space
$(\mathcal{M}_{0}^{(n)},\goth{B}(\mathcal{M}_{0}^{(n)}))$. We remind
that if the intensity measures $\Lambda_{\vartheta_1 }^{(n)}$ and
$\Lambda_{\vartheta_2}^{(n)}$ are equivalent then the corresponding
probability measures $\bP_{\vartheta_1}^{(n)}$ and
$\bP_{\vartheta_2}^{(n)}$ are equivalent and the likelihood ratio is
given by
\begin{align*}
&\frac{{\rm d}\bP^{(n)}_{\vartheta_{2}}}{{\rm
d}\bP^{(n)}_{\vartheta_{1}}}\left(X^{\left(n\right)}\right)=\\&=
\exp\left\{\int_{\mathbb{A}_n}\ln\frac{S(\vartheta_2,x)}
{S(\vartheta_1,x)}\;X^{(n)}({\rm
d}x)-\int_{\mathbb{A}_n}\left[S(\vartheta_2,x)-S(\vartheta_1,x)\right]\,
{\rm d}x\right\}.
\end{align*}
For a proof see \cite{Kut98} page 28. Here the stochastic integral
$$
\int_{\mathbb{A}_n}f(x)\;X^{(n)}({\rm d}x)=\sum_{x_i\in
\mathbb{A}_n} f(x_i)
$$
where $\left\{x_i\right\}$ are the {\sl events} (random points) of
the Poisson process.

\bigskip

Since the main tool used in this work is based on the Edgeworth
expansion, here we present the conditions under which the
distribution function
\begin{align*}
F_{n}(y)=\bP_{\vartheta}^{\left(n\right)}\left\{ I_n(f)<y\right\},
\end{align*}
of the stochastic integral
$$
I_n(f)=\int_{\mathbb{A}_{n}} f_n\left(x\right)\;\pi^{(n)}({\rm d}x),
$$
admits an Edgeworth type expansion for two terms after the Gaussian
term, where $\pi ^{(n)}\left({\rm d}x\right)=X^{(n)}\left({\rm
d}x\right)-S(\vartheta,x)\;{\rm d}x$ is the centered Poisson
process. We suppose that
\begin{align*}
\label{26} \int_{\mathbb{A}_{n}}f_{n}(x) ^{2}\; S(\vartheta,x)
\;{\rm d}x=1.
\end{align*}
The expansion is obtained under the following two conditions:
\begin{itemize}
\item[${\cal B}_1.$] {\sl There exists a sequence of real numbers
$\varepsilon_{n}\rightarrow 0$, as $n\rightarrow\infty$ and
constants $C_{r}>0, r=3,4,5$,  such that}
$$
\int_{\mathbb{A}_{n}} \abs{f_{n}(x)}^{r} \; S(\vartheta,x)\;{\rm
d}x\leq C_{r}\;\varepsilon_{n}^{r-2}.
$$
\item[${\cal B}_2.$] {\sl There exist constants  $\gamma\geq 5/2$ and
$c_{0}>0$  satisfying the inequality\;
$\frac{C_{3}}{3!}c_{0}+\frac{C_{4}}{4!}c_{0}^{2}+\frac{C_{5}}{5!}c_{0}^{3}-\frac{1}{2}<0$
such that
$$
\inf_{\frac{c_{0}\varepsilon_{n}^{-1}}{2}<t<\frac{\varepsilon_n^{-2}}{2}}
\int_{\mathbb{A}_{n}}\sin^{2}\left(t
f_{n}(x)\right)\;S(\vartheta,x)\;{\rm d}x\geq\gamma
\ln\varepsilon_{n}^{-1}
$$
for all large $n$.}
\end{itemize}
Let us introduce the cumulants:
$$
\gamma_{r,n}=\int_{\mathbb{A}_{n}} f_n(x)^r\;S(\vartheta,x)\;{\rm
d}x,\qquad r=3,4
$$
and the Hermit polynomials:
$$
H_2(y)=y^2-1,\quad H_3(y)=y^3-3y,\quad H_5(y)=y^5-10y^3+15y.
$$
\begin{theorem}
\label{EX} Let the conditions ${\cal B}_1,\;{\cal B}_2$ be
fulfilled, then uniformly in $y\in \mathbb{R}$
\begin{align*}
F_{n}(y)={\cal N}(y)-\frac{\gamma_{3,n}}{3!}\,H_2(y)\,n(y)
-\frac{\gamma_{4,n}}{4!}\,H_3(y)\,n(y)-\frac{\gamma_{3,n}^2}{72}\,H_5(y)\,n(y)
+O(\varepsilon_{n}^{3}),
\end{align*}
for all $n$ large. Here ${\cal N}(y)$ and $n(y)$ denote the
distribution and density functions of the standard Gaussian law, respectively.
\end{theorem}
\noindent Note that $\gamma_{r,n}=O(\varepsilon_{n}^{r-2}),\; r=3,4$
by ${\cal B}_1$.

\noindent {\bf Proof.} See \cite{Fa}, Page 36. The proof is a
special case of a general theorem given by Kutoyants, where the
expansion is obtained by the powers of $\varepsilon_n$ up to order
$\varepsilon_n^k$, $k=1,2,...$ (see (\cite{Kut98}), page 131).

\begin{corollary}
\label{tt} Let $0<\alpha<1$ be given and the conditions
$\mathcal{B}_1-\mathcal{B}_2$ be fulfilled. Then the equation
$F_n(y)=1-\alpha+O(\varepsilon_n^3)$ has a solution
$y=c_{n,\alpha}$,
$$
c_{n,\alpha}=z_\alpha+\frac{\gamma_{3,n}}{6}\,H_2(z_\alpha)
+\frac{\gamma_{4,n}}{24}\,H_3(z_\alpha)+\frac{\gamma_{3,n}^2}{72}\,H_5(z_\alpha)
.
$$
\end{corollary}
\noindent See \cite{Fa}, page 40 for proof.

\section{Second Order Efficiency}
\noindent Let $\beta_{n}\left(u, \phi_{n}\right)$ denote the power
of a test $\phi_n$ at the local alternative
$\vartheta_u=\vartheta_0+\varphi_n u$, {\sl i.e.,}
$$
\beta_{n}\left(u, \phi_{n}\right)=\bE _{\vartheta_u}\phi_n(X^{(n)}).
$$
Under certain regularity conditions, slightly weaker than
$\mathcal{D}_1-\mathcal{D}_3$ (see the next section and \cite{Fa1}),
the score test
\begin{equation*}
\phi_n^*\left(X^{\left(n\right)}\right)=\begin{cases}1,& \,
{\rm if}\,\, \Delta_{n}\left(\vartheta_{0}\right)>c_n\\
0,& \, {\rm if}\,\, \Delta_{n}\left(\vartheta_{0}\right)\leq c_n,
\end{cases}
\end{equation*}
based on
$$
\Delta_{n}\left(\vartheta_{0}\right)=\varphi_{n}
\int_{\mathbb{A}_{n}}\frac{\dot{S}\left(\vartheta_{0}, x
\right)}{S\left(\vartheta_{0}, x \right)}\,\pi^{(n)}\left({\rm
d}x\right)
$$
with
\begin{equation*}
 c_n=z_{\alpha}-
\frac{\gamma_{3,n}}{6}\left(1-z_{\alpha}^{2}\right),\quad
\gamma_{3,n}=\varphi_{n}^3\int_{\mathbb{A}_{n}}
\frac{\dot{S}\left(\vartheta_{0}, x \right)^{3}}
{S\left(\vartheta_{0}, x \right)^{2}}\;{\rm d}x
\end{equation*}
is second order efficient, {\sl i.e.,} it satisfies
$$
\sup_{0\leq u \leq K}\abs{\beta_{n}\left(u,
\phi_{n}^{*}\right)-\beta_{n}(u,
\tilde{\phi}_{n})}=O(\varepsilon_n^2),
$$
for any $K>0$. For proof see \cite{Fa1}, Theorem 6. Indeed to
establish the second order efficiency we obtain the following
representation of distribution functions of
$\Delta_{n}\left(\vartheta_{0}\right)$ and $\Lambda_n(u)$ under the
local alternative $\vartheta_u$:
\begin{align*}
\beta_{n}\left(u, \phi_{n}^{*}\right)&=
\mathcal{N}\left(\frac{m_n(u)-c_n}{\eta_n}\right)-
\frac{\gamma_{3,n}(u)}{6}\,(1-(u-z_{\alpha})^2)\,n(u-z_{\alpha})+
O(\varepsilon_n^2)\\
\nonumber\beta_{n}(u,
\tilde{\phi}_{n})&=\mathcal{N}\left(\frac{\mu_n(u)-b_n(u)}{\sigma_n(u)}\right)-
\frac{\gamma'_{3,n}(u)}{6}(1-(u-z_{\alpha})^2)n(u-z_{\alpha})+
O(\varepsilon_n^2),
\end{align*}
where
\begin{align*}
m_n(u)&=\Ex_{\vartheta_u}\Delta_{n}\left(\vartheta_{0}\right)=
\varphi_{n}\int_{\mathbb{A}_{n}} \frac{\dot{S}\left(\vartheta_{0}, x
\right)} {S\left(\vartheta_{0}, x
\right)}\left(S\left(\vartheta_{u}, x
\right)-S\left(\vartheta_{0}, x \right)\right)\;{\rm d}x,\\
\eta_n^2&=\Ex_{\vartheta_u}\left(\Delta_{n}(\vartheta_{0})-m_n(u)\right)^2=
\varphi_{n}^2\int_{\mathbb{A}_{n}} \frac{\dot{S}\left(\vartheta_{0},
x \right)^2} {S\left(\vartheta_{0}, x
\right)^2}S\left(\vartheta_{u}, x
\right)\;{\rm d}x,\\
\mu_n(u)&=\Ex_{\vartheta_u}\Lambda_n(u)=\int_{\mathbb{A}_{n}}\left(
\ln\frac{S\left(\vartheta_{u},x\right)}{S\left(\vartheta_{0}, x
\right)}S\left(\vartheta_{u},x \right)-S\left(\vartheta_{u}, x
\right)+S\left(\vartheta_{0},x\right)\right){\rm d}x,\\
\sigma^2_n(u)&=\Ex_{\vartheta_u}\left(\Lambda_n(u)-\mu_n(u)\right)^2=\int_{\mathbb{A}_{n}}
\left(\ln\frac{S\left(\vartheta_{u}, x
\right)}{S\left(\vartheta_{0},
x \right)}\right)^2\,S\left(\vartheta_{u}, x \right){\rm d}x,\\
\gamma_{3,n}(u)&=\frac{\varphi_{n}^3}{\eta_n^3}\int_{\mathbb{A}_{n}}
\frac{\dot{S}\left(\vartheta_{0}, x \right)^{3}}
{S\left(\vartheta_{0}, x \right)^{3}}S\left(\vartheta_{u}, x
\right)\;{\rm
d}x,\\
\gamma'_{3,n}(u)&=\frac{1}{\sigma_n(u)^3}\int_{\mathbb{A}_{n}}
\left(\ln\frac{S\left(\vartheta_{u}, x
\right)}{S\left(\vartheta_{0}, x
\right)}\right)^3\,S\left(\vartheta_{u}, x \right){\rm d}x.
\end{align*}
Hence it suffices to show that
\begin{align*}
\sup_{0< u \leq K}\abs{\gamma_{3,n}(u)-\gamma_{3,n}'(u)}&=O(\varepsilon_n^2),\\
\sup_{0< u \leq K}\abs{\frac{m_n(u)-c_n}{\eta_n}-
\frac{\mu_n(u)-b_n(u)}{\sigma_n(u)}}&=O(\varepsilon_n^2).
\end{align*}
The case $u=0$ corresponds to the size of  $\tilde{\phi}_{n}$ which
is equal to $\alpha$ and the size of $\phi_{n}^{*}$ given by
$\alpha+O(\varepsilon_n^2)$. The latter follows from the Edgeworth
expansion of distribution function of
$\Delta_{n}\left(\vartheta_{0}\right)$ under the null hypothesis.

\bigskip

$\mathbf{Representation\; of\; the\; power\;}$. Now we obtain the
explicit representation of the power of $\phi_{n}^{*}$ (up to order
$O(\varepsilon_n^2)$). The Taylor expansion
$$
S(\vartheta_{u},x)=S(\vartheta_{0},x)+ \varphi_n
u\dot{S}(\vartheta_{0},x)+\frac{\varphi_n^2\,u^2}{2}\ddot{S}(\vartheta_{0},x)
+\frac{\varphi_n^3\,u^3}{3!} S^{(3)}(\vartheta_{n},x),
$$
for some intermediate point
$\vartheta_0<\vartheta_n=\vartheta_n(u,x)<\vartheta_{u}$, implies
that
$$
m_n(u)-c_n=u-z_{\alpha}+\frac{\varphi_n^3\,u^2}{2}
\int_{\mathbb{A}_n}\frac{\dot{S}(\vartheta_{0},x)\,\ddot{S}(\vartheta_{0},x)}
{S(\vartheta_{0},x)}\;{\rm
d}x+\frac{\gamma_{3,n}}{6}(1-z_{\alpha}^2)+O(\varepsilon_n^2).
$$
Similarly we obtain
$$
\eta^2_n=1+\varphi_n^3
u\int_{\mathbb{A}_n}\frac{\dot{S}(\vartheta_{0},x)^3}
{S(\vartheta_{0},x)^2}\;{\rm
d}x+O(\varepsilon_n^2)=1+u\,\gamma_{3,n}+O(\varepsilon_n^2).
$$
Hence
\begin{align*}
\frac{m_n(u)-c_n}{\eta_n}=u-z_{\alpha}
&+\frac{1-z_{\alpha}^2-3u\,(u-z_{\alpha})}{6}\;\gamma_{3,n}
+\\&+\frac{\varphi_n^3\,u^2}{2}
\int_{\mathbb{A}_n}\frac{\dot{S}(\vartheta_{0},x)\,\ddot{S}(\vartheta_{0},x)}
{S(\vartheta_{0},x)}\;{\rm d}x+O(\varepsilon_n^2).
\end{align*}
On the other hand by the help of Taylor expansion we obtain
$$
\gamma_{3,n}(u)=\gamma_{3,n}+O(\varepsilon_n^2).
$$
Therefore we obtain the following representation
\begin{equation}
\label{power21} \beta_{n}\left(u,
\phi_{n}^{*}\right)=\mathcal{N}\left(u-z_{\alpha}\right)
+Q_n(u)\,n(u-z_{\alpha})+O(\varepsilon_n^2),
\end{equation}
for any $u>0$, where the polynomial (in $u$)
$$
Q_n(u)=\frac{u(z_{\alpha}-2u)}{6}\;\gamma_{3,n}+\frac{\varphi_n^3\,u^2}{2}
\int_{\mathbb{A}_n}\frac{\dot{S}(\vartheta_{0},x)\,\ddot{S}(\vartheta_{0},x)}
{S(\vartheta_{0},x)}\;{\rm d}x,
$$
is of order $O(\varepsilon_n)$. The equation \eqref{power21},
refines the first order representation \eqref{power1}. Notice also
that the the second order efficiency of $\phi_{n}^{*}$ implies that
we have
the same representation for the power of $\tilde{\phi}_n$.%

\section{Power Loss}
The power function of a second order efficient test agrees with that
of the most powerful test up to order $O(\varepsilon_{n}^2)$. Hence
it is natural to consider the {\sl power loss} of $\phi_{n}^{*}$,
which is defined for any $u>0$ by
\begin{equation*}
r(u)=\lim_{n\rightarrow\infty}
\varepsilon_{n}^{-2}\left(\beta_{n}(u, \tilde{\phi}_n)-
\beta_{n}\left(u, \phi_{n}^{*}\right)\right).
\end{equation*}
\noindent Below $S^{(j)}(\vartheta,x)$ denotes the $j$th derivative
of $S(\vartheta,x)$ with respect to $\vartheta$. We write
$\dot{S}(\vartheta,x)$ and $\ddot{S}(\vartheta,x)$ for the first and
second derivatives, respectively. We consider the following
conditions:
\begin{itemize}
\item[$\mathcal{D}_1.$]  The intensity function $S(\vartheta,x)$ is four times differentiable with respect to
$\vartheta$ in a right neighborhood of $\vartheta_0$.
\item[$\mathcal{D}_2.$] The conditions
$\mathcal{B}_1$ and $\mathcal{B}_2$ are satisfied for the stochastic
integrals $\Delta_n(\vartheta_0)$ and $\Lambda_n(u)$ under ${\cal
H}_{0}$ and ${\cal H}_{u}$ with some sequence
$\varepsilon_n\rightarrow 0$,
\item[$\mathcal{D}_3.$] There exists some functions
$f_j(x), j=0,...,4,\;x\in \mathbb{A}_n$ not depending on $\vartheta$
such that $S(\vartheta,x)\geq f_0(x), \abs{S^{(j)}(\vartheta,x)}\leq
f_j(x), j=0,...,4$ for all $x\in\mathbb{A}_n$ and all $\vartheta$ in
a right neighborhood of $\vartheta_0$. We suppose also that
\begin{align*}
\varphi_n^k\int_{\mathbb{A}_n}\frac{\abs{f_1(x)}^{k}}
{f_0(x)^{k-1}}\;{\rm d}x &=O(\varepsilon_n^{k-2}),\quad k=2,3,4\\
\varphi_n^{2j}\int_{\mathbb{A}_n}\frac{f_j(x)^{2}} {f_0(x)}\;{\rm
d}x &=O(\varepsilon_n^{2j-2}),\quad j=2,3,4.
\end{align*}
\end{itemize}
{\bf Example 1.}  Let $X^{(n)}$ be a realization of a Poisson
process on the set $\mathbb{A}_n=[\,0,n]$ with positive intensity
function $S(\vartheta,x)=\vartheta S(x)+\lambda$ (amplitude
parameter) or $S(\vartheta,x)= S(\vartheta+x)+\lambda$ (phase
parameter), where $S(\cdot)$ is a two times differentiable periodic
function and $\lambda>0$ (dark current) is a known constant. In both
cases the condition $\mathcal{D}_3$ is satisfied with $\varphi_n\sim
C\,n^{-1/2}$ for some $C>0$ and $\varepsilon_n= n^{-1/2}$. For the
frequency modulation model $S(\vartheta,x)=S(\vartheta\, x)+\lambda$
we have $\varphi_n\sim C\,n^{-3/2}$ and $\varepsilon_n= n^{-1/2}$.

Let $J_n$ denote the quantity:
\begin{align*}
J_n&= \varphi_n^4\int_{\mathbb{A}_{n}}
\frac{\left(\dot{S}\left(\vartheta_{0},
x\right)^{2}-S\left(\vartheta_{0},
x\right)\ddot{S}\left(\vartheta_{0},
x\right)\right)^{2}}{S\left(\vartheta_{0}, x \right)^{3}}{\rm
d}x-\\&-\left(\varphi_n^3\int_{\mathbb{A}_{n}}
\frac{\dot{S}\left(\vartheta_{0},
x\right)\left(\dot{S}\left(\vartheta_{0},
x\right)^{2}-S\left(\vartheta_{0},
x\right)\ddot{S}\left(\vartheta_{0},
x\right)\right)}{S\left(\vartheta_{0}, x \right)^{3}}{\rm
d}x\right)^{2}.
\end{align*}
Note that $J_n=O(\varepsilon_n^2)$, by $\mathcal{D}_3$. We have the
following theorem:
\begin{theorem}
Let the conditions $\mathcal{D}_1-\mathcal{D}_3$ be fulfilled. Then
the power loss of $\phi^*_n$ with respect to the most powerful test
$\tilde{\phi}_n$ is equal to
\begin{equation*}
r(u)=\frac{u^3\,n(u-z_\alpha)}{8}
\,\lim_{n\rightarrow\infty}\left(\varepsilon_{n}^{-2}\,J_n\right),
\end{equation*}
for any $u>0$.
\end{theorem}
{\bf Proof.} By $\mathcal{D}_2$ we can write the following third
order expansions:
\begin{align*}
\beta_{n}\left(u, \phi_{n}^{*}\right)&= {\cal
N}(a_n)+\frac{\gamma_{3,n}\left(u\right)}{3!}\,H_2(a_n)\,n(a_n)
-\frac{\gamma_{4,n}(u)}{4!}\,H_3(a_n)\,n(a_n)-\\&-\frac{\gamma^2_{3,n}(u)}{72}\,H_5(a_n)\,n(a_n)
+O(\varepsilon_{n}^{3})\\
\beta_{n}(u, \tilde{\phi}_n)&={\cal
N}(A_n)+\frac{\gamma'_{3,n}\left(u\right)}{3!}\,H_2(A_n)\,n(A_n)
-\frac{\gamma'_{4,n}(u)}{4!}\,H_3(A_n)\,n(A_n)-\\&-\frac{\gamma'_{3,n}(u)^2}{72}\,H_5(A_n)\,n(A_n)
+O(\varepsilon_{n}^{3}),
\end{align*}
where
\begin{align*}
a_n&=\frac{m_n(u)-c_n}{\eta_n},\quad\qquad A_n=\frac{\mu_n(u)-b_n(u)}{\sigma_n(u)}\\
\gamma_{r,n}(u)&=\frac{\varphi_{n}^r}{\eta_n^r}\int_{\mathbb{A}_{n}}
\frac{\dot{S}\left(\vartheta_{0}, x \right)^{r}}
{S\left(\vartheta_{0}, x \right)^{r}}S\left(\vartheta_{u}, x
\right)\;{\rm d}x,\qquad r=3,4\\
\gamma'_{r,n}(u)&=\frac{1}{\sigma_n(u)^r}\int_{\mathbb{A}_{n}}
\left(\ln\frac{S\left(\vartheta_{u}, x
\right)}{S\left(\vartheta_{0}, x
\right)}\right)^r\,S\left(\vartheta_{u}, x \right){\rm d}x, \qquad
r=3,4.
\end{align*}
Using the Taylor expansions of $S\left(\vartheta_{u}, x \right)$ and
$\ln\frac{S\left(\vartheta_{u}, x \right)}{S\left(\vartheta_{0}, x
\right)}$, we get
\begin{align*}
\gamma_{4,n}(u)&=\gamma_{4,n}+O(\varepsilon_{n}^{3}), \quad
\gamma_{3,n}(u)^2=\gamma_{3,n}^2+O(\varepsilon_{n}^{3})\\
\gamma'_{4,n}(u)&=\gamma_{4,n}+O(\varepsilon_{n}^{3}), \quad
\gamma'_{3,n}(u)^2=\gamma_{3,n}^2+O(\varepsilon_{n}^{3}).
\end{align*}
Since $A_n-a_n=O(\varepsilon_{n}^{2})$ (which follows from the
second order efficiency of $\phi_n^*$), then
$$
\mathcal{N}\left(A_n\right)=\mathcal{N}\left(a_n\right)+\left(A_n-a_n\right)\,n(a_n)+O(\varepsilon_{n}^{4}).
$$
Therefore
\begin{align*}
\beta_{n}\left(u, \phi_{n}^{*}\right)&= {\cal
N}(a_n)+\frac{\gamma_{3,n}}{3!}\,H_2(a_n)\,n(a_n)
-\frac{\gamma_{4,n}}{4!}\,H_3(a_n)\,n(a_n)-\\&-\frac{\gamma^2_{3,n}}{72}\,H_5(a_n)\,n(a_n)
+O(\varepsilon_{n}^{3}),\\
\beta_{n}(u,
\tilde{\phi}_n)&=\mathcal{N}\left(a_n\right)+\left(A_n-a_n\right)\,n(a_n)
+\frac{\gamma'_{3,n}\left(u\right)}{3!}\,H_2(a_n)\,n(a_n)
-\\&-\frac{\gamma_{4,n}}{4!}\,H_3(a_n)\,n(a_n)-\frac{\gamma_{3,n}^2}{72}\,H_5(a_n)\,n(a_n)
+O(\varepsilon_{n}^{3}).
\end{align*}
Letting $\Delta=\Delta(u)=u-z_\alpha$ and taking into account the
fact that $a_n-\Delta=O(\varepsilon_{n})$, we can write
\begin{align*}
\beta_{n}(u, \tilde{\phi}_n)-\beta_{n}\left(u,
\phi_{n}^{*}\right)&=\left(A_n-a_n\right)n(\Delta)+\\&+
\frac{\,H_2(\Delta)\,n(\Delta)}{6}
\left(\gamma'_{3,n}\left(u\right)-\gamma_{3,n}\left(u\right)\right)
+O(\varepsilon_{n}^{3}).
\end{align*}
Therefore we have to consider the terms $A_n-a_n$ and
$\gamma'_{3,n}\left(u\right)-\gamma_{3,n}\left(u\right)$. Since the
investigation is at accuracy level $O(\varepsilon_{n}^{3})$, we
modify the threshold $c_n$ as follows:
$$
c_n=z_{\alpha}+\frac{\gamma_{3,n}}{6}\,H_2(z_{\alpha})+
\frac{\gamma_{4,n}}{4!}\,H_3(z_{\alpha})+
\frac{\gamma_{3,n}^2}{72}\,H_5(z_{\alpha}).
$$
With this constant the probability of error of the first kind of
$\phi^*_n$ is equal to $\alpha+O(\varepsilon_n^3)$. Similarly the
explicit form of  $b_n(u)$, the threshold of $\tilde{\phi}_n$, can
be written as:
\begin{align*}
b_n(u)=\mu_n+\sigma_n\left(z_{\alpha}+
\frac{\gamma'_{3,n}}{6}\,H_2(z_{\alpha})+
\frac{\gamma'_{4,n}}{4!}\,H_3(z_{\alpha})+ \frac{\gamma'^{2}
_{3,n}}{72}\,H_5(z_{\alpha})\right)+O(\varepsilon_{n}^{3}),
\end{align*}
where
\begin{align*}
\mu_n&=\Ex_{\vartheta_0}\Lambda_n(u)=\int_{\mathbb{A}_{n}}\left[
\ln\frac{S\left(\vartheta_{u},x\right)}{S\left(\vartheta_{0}, x
\right)}-\frac{S\left(\vartheta_{u}, x
\right)}{S\left(\vartheta_{0},x
\right)}+1\right]S\left(\vartheta_{0},x
\right){\rm d}x,\\
\sigma_n^2&=\Ex_{\vartheta_0}\left(\Lambda_n(u)-\mu_n\right)^2=\int_{\mathbb{A}_{n}}
\left(\ln\frac{S\left(\vartheta_{u}, x
\right)}{S\left(\vartheta_{0},
x \right)}\right)^2\,S\left(\vartheta_{0}, x \right){\rm d}x,\\
\gamma'_{r,n}&=\frac{1}{\sigma_n^r}\int_{\mathbb{A}_{n}}
\left(\ln\frac{S\left(\vartheta_{u}, x
\right)}{S\left(\vartheta_{0}, x
\right)}\right)^r{S\left(\vartheta_{0}, x \right)}\;{\rm d}x, \quad
r=3,4.
\end{align*}
Below for the sake of simplicity we use the following notations:
$$
I(r_0,r_1,r_2)=\int_{\mathbb{A}_{n}}
\frac{\dot{S}\left(\vartheta_{0},
x\right)^{r_1}\ddot{S}\left(\vartheta_{0},
x\right)^{r_2}}{S\left(\vartheta_{0}, x \right)^{r_0}}{\rm d}x
$$
for nonnegative integers $r_0,r_1,r_2$. The long but straightforward
calculations show that (see Appendix)
\begin{align}
\label{A-a} A_n-a_n&=\frac{u^3}{8}\,\varphi_n^4\,I(1,0,2)\nonumber
-\frac{2u^3-2u^2z_{\alpha}-u(1-z_{\alpha}^2)}{4}\,\varphi_n^4\,I(2,2,1)+\\&+
\frac{9u^3-12u^2z_{\alpha}-6u(1-z_{\alpha}^2)}{24}\,\varphi_n^4\,I(3,4,0)+\\&+\nonumber
\frac{9u^3-6u^2z_{\alpha}+2u(1-z_{\alpha}^2)}{24}\,\varphi_n^6\,\,I^{2}(2,3,0)+\\&+\nonumber
\frac{6u^3-6u^2z_{\alpha}-3u(1-z_{\alpha}^2)}{12}\,\varphi_n^6\,I(2,3,0)\,I(1,1,1)-\\&-\nonumber
\frac{u^3}{8}\,\varphi_n^6\,I(1,1,1)^{2}+O(\varepsilon_n^3).
\end{align}
Also for the term
$\gamma'_{3,n}\left(u\right)-\gamma_{3,n}\left(u\right)$ we find
that
\begin{align}
\label{gamma}\gamma'_{3,n}\left(u\right)-\gamma_{3,n}\left(u\right)&=\frac{3u}{2}\,\varphi_n^4\,I(2,2,1)+
\frac{3u}{2}\,\varphi_n^6\,I(2,3,0)^{2}-\frac{3u}{2}\,\varphi_n^4\,I(3,4,0)\nonumber-\\&-
\frac{3u}{2}\,\varphi_n^6\,I(2,3,0)\,I(1,1,1)+O(\varepsilon_n^3).
\end{align}
See Appendix for some details. Combining these results yield
\begin{align*}
&n^{-1}(\Delta)\left(\beta_{n}(u, \tilde{\phi}_n)-\beta_{n}\left(u,
\phi_{n}^{*}\right)\right)=\frac{u^3}{8}\,\varphi_n^4\,I(1,0,2)-\frac{u^3}{4}\,\varphi_n^4\,I(2,2,1)
+\\&+
\frac{u^3}{8}\,\varphi_n^4\,I(3,4,0)-\frac{u^3}{8}\,\varphi_n^6\,I(2,3,0)^{2}
-\frac{u^3}{8}\,\varphi_n^6\,I(1,1,1)^{2}+\\&
+\frac{u^3}{4}\,\varphi_n^6\,I(2,3,0)\,I(1,1,1)+O(\varepsilon_{n}^{3})=
\frac{u^3}{8}\,\varphi_n^4\left[\,I(1,0,2)-2\,I(2,2,1)+\right.\\&+
\left.I(3,4,0)\right]-\frac{u^3}{8}
\,\varphi_n^6\left[I(2,3,0)-I(1,1,1)\right]^2+O(\varepsilon_{n}^{3}).
\end{align*}
Hence for any $u>0$
\begin{align*}
\beta_{n}(u, \tilde{\phi}_n)-\beta_{n}\left(u,
\phi_{n}^{*}\right)=\frac{u^3\,n(\Delta)}{8}\,J_n+O(\varepsilon_{n}^{3})
\end{align*}
where
\begin{align*}
J_n&=\,\varphi_n^4\left[\,I(1,0,2)-2\,I(2,2,1)+
I(3,4,0)\right]-\varphi_n^6\,\left[I(2,3,0)-I(1,1,1)\right]^2.
\end{align*}
This can be simplified as
\begin{align*}
J_n&= \varphi_n^4\int_{\mathbb{A}_{n}}
\frac{\left(\dot{S}\left(\vartheta_{0},
x\right)^{2}-S\left(\vartheta_{0},
x\right)\ddot{S}\left(\vartheta_{0},
x\right)\right)^{2}}{S\left(\vartheta_{0}, x \right)^{3}}{\rm
d}x-\\&-\left(\varphi_n^3\int_{\mathbb{A}_{n}}
\frac{\dot{S}\left(\vartheta_{0},
x\right)\left(\dot{S}\left(\vartheta_{0},
x\right)^{2}-S\left(\vartheta_{0},
x\right)\ddot{S}\left(\vartheta_{0},
x\right)\right)}{S\left(\vartheta_{0}, x \right)^{3}}{\rm
d}x\right)^{2}.
\end{align*}
Note that by the Cauchy-Schwartz inequality;
\begin{align*}
&\left(\varphi_n^3\int_{\mathbb{A}_{n}}
\frac{\dot{S}\left(\vartheta_{0},
x\right)\left(\dot{S}\left(\vartheta_{0},
x\right)^{2}-S\left(\vartheta_{0},
x\right)\ddot{S}\left(\vartheta_{0},
x\right)\right)}{S\left(\vartheta_{0}, x \right)^{3}}{\rm
d}x\right)^{2}\leq \\\leq &\left(\varphi_n^2\int_{\mathbb{A}_{n}}
\frac{\dot{S}\left(\vartheta_{0}, x\right)^2}{S\left(\vartheta_{0},
x \right)}{\rm d}x\right)\left(\varphi_n^4\int_{\mathbb{A}_{n}}
\frac{\left(\dot{S}\left(\vartheta_{0},
x\right)^{2}-S\left(\vartheta_{0},
x\right)\ddot{S}\left(\vartheta_{0},
x\right)\right)^{2}}{S\left(\vartheta_{0}, x \right)^{3}}{\rm
d}x\right)=\\&=\varphi_n^4\int_{\mathbb{A}_{n}}
\frac{\left(\dot{S}\left(\vartheta_{0},
x\right)^{2}-S\left(\vartheta_{0},
x\right)\ddot{S}\left(\vartheta_{0},
x\right)\right)^{2}}{S\left(\vartheta_{0}, x \right)^{3}}{\rm d}x.
\end{align*}
Hence as one expect $J_n\geq 0$, because the test $\tilde{\phi}_n$
has the maximum of power at the local alternative $\vartheta_u$.
Furthermore $\mathcal{D}_3$ implies that the quantity $J_n$ is of
order $O(\varepsilon_{n}^{2})$. Therefore
$$
\varepsilon_{n}^{-2}\left(\beta_{n}(u, \tilde{\phi}_n)-
\beta_{n}\left(u,
\phi_{n}^{*}\right)\right)=\frac{u^3\,n(\Delta)}{8}
\,\varepsilon_{n}^{-2}\,J_n+O(\varepsilon_n)
$$
which completes the proof of the theorem.

\bigskip

{\bf Example 2. (Amplitude Parameter.)} Suppose that we observe a
realization $X^{(n)}$ of a Poisson process on the set
$\mathbb{A}_{n}=[\,0, n]$, $n=1,2,\cdots$ with the intensity
function
$$
S(\vartheta,x)=\vartheta\,S(x)+\lambda, \qquad\vartheta>0
$$
where $\lambda$ is a known positive constant (dark-current) and
$S(x)$ is a known, nonconstant, differentiable (with respect to $x$)
and periodic function with period $\tau>0$. We have two hypotheses
${\cal H}_{0}:\vartheta=\vartheta_{0}$ against ${\cal H}_{1}:
\vartheta>\vartheta_{0}$, where $\vartheta_0>0$. The intensity
function $S(\vartheta,x)$ is supposed to be positive in a right
neighborhood of $\vartheta_0$ and all $x$. We obtain the power loss
of the test
$$
\phi^{*}_{n}\left(X^{(n)}\right)=\chi_{\left\{\Delta_{n}\left(\vartheta_{0}\right)>
c_{n}\right\}}
$$
with
$$ c_n=z_\alpha+\frac{\gamma_{3,n}}{6}\,H_2(z_\alpha)
+\frac{\gamma_{4,n}}{24}\,H_3(z_\alpha)+\frac{\gamma_{3,n}^2}{72}\,H_5(z_\alpha),
$$
where
\begin{align*}
\gamma_{3,n}&=\frac{1}{\tau\,A^3\,\sqrt{n}}\int_{0}^{\tau}\frac{S(x)^3}{\left(\vartheta_0
S(x)+\lambda\right)^2}\;{\rm d}x\\
\gamma_{4,n}&=\frac{1}{\tau\,A^{4}\,n}\int_{0}^{\tau}\frac{S(x)^4}{\left(\vartheta_0
S(x)+\lambda\right)^3}\;{\rm d}x\\
A^{2}&=\frac{1}{\tau}\int_{0}^{\tau}\frac{S(x)^2}{\vartheta_0
S(x)+\lambda}\;{\rm d}x.
\end{align*}
The conditions $\mathcal{D}_1-\mathcal{D}_3$ are satisfied (see
\cite{Fa1}, pp. 202-204) with $\varepsilon_n=n^{-1/2}$ and
$\varphi_n=n^{-1/2}\,A^{-1}(1+O(n^{-1}))\rightarrow 0$. One can
obtain
$$
\lim_{n\rightarrow\infty}\left(\varepsilon_{n}^{-2}\,J_n\right)
=\tau\left(\frac{\int_{0}^{\tau}\frac{S(x)^4}{\left(\vartheta_0
S(x)+\lambda\right)^3}\;{\rm
d}x}{\left(\int_{0}^{\tau}\frac{S(x)^2}{\vartheta_0
S(x)+\lambda}\;{\rm
d}x\right)^2}-\frac{\left(\int_{0}^{\tau}\frac{S(x)^3}{\left(\vartheta_0
S(x)+\lambda\right)^3}\;{\rm
d}x\right)^2}{\left(\int_{0}^{\tau}\frac{S(x)^2}{\vartheta_0
S(x)+\lambda}\;{\rm d}x\right)^3}\right)\equiv B.
$$
Hence the power loss is equal to
$$
r(u)=\frac{u^3\,n(u-z_{\alpha})}{8}\,B.
$$

\noindent {\bf Example 3. (Frequency parameter)} In this example we
consider a strongly nonhomogeneous case with nonclassical rate
$n^{-3/2}$ (instead of $n^{-1/2}$ in the i.i.d. case). Suppose that
we observe a realization $X^{(n)}$ of a Poisson process on the set
$\mathbb{A}_{n}=[\,0, n]$, $n=1,2,\cdots$ with periodic intensity
function
$$
S(\vartheta,x)=e^{\sin(\vartheta \,x)}, \qquad\vartheta>0.
$$
We have the simple hypothesis ${\cal H}_{0}:\vartheta=\vartheta_{0}$
against ${\cal H}_{1}: \vartheta>\vartheta_{0}$. We consider the
test $
\phi^{*}_{n}\left(X^{(n)}\right)=\chi_{\left\{\Delta_{n}\left(\vartheta_{0}\right)>
c_{n}\right\}} $ based on the statistic
$$
\Delta_{n}\left(\vartheta_{0}\right)=\varphi_n\int_{0}^{n}x\,\cos(\vartheta_0\,x)\left(X^{(n)}\left({\rm
d}x\right)-e^{\sin(\vartheta_0 \,x)}{\rm d}x\right)
$$
and the threshold
$$ c_n=z_\alpha+\frac{\gamma_{3,n}}{6}\,H_2(z_\alpha)
+\frac{\gamma_{4,n}}{24}\,H_3(z_\alpha)+\frac{\gamma_{3,n}^2}{72}\,H_5(z_\alpha),
$$
where letting $ \tau=\frac{2\pi}{\vartheta_0}$ we have
\begin{align*}
\gamma_{3,n}&=\frac{C^3}{4\tau\,\sqrt{n}}\int_{0}^{\tau}\cos^3(\vartheta_0\,x)\,e^{\sin(\vartheta_0
\,x)}\;{\rm d}x\\
\gamma_{4,n}&=\frac{C^4}{5\tau\,n}\int_{0}^{\tau}\cos^4(\vartheta_0\,x)\,e^{\sin(\vartheta_0
\,x)}\;{\rm d}x\\
C^{-2}&=\frac{1}{3\,\tau}\int_{0}^{\tau}\cos^2(\vartheta_0\,x)\,e^{\sin(\vartheta_0
\,x)}\;{\rm d}x.
\end{align*}
The conditions $\mathcal{D}_1-\mathcal{D}_3$ are satisfied (see
\cite{Fa1}, pp. 205-207) with $\varepsilon_n=n^{-1/2}$ and
$\varphi_n\sim C\,n^{-3/2}$. 

The power loss is equal to
\begin{equation*}
r(u)=\frac{9\,\tau\,u^3\,n(u-z_\alpha)}{40}
\frac{\int_{0}^{\tau}\sin^2(\vartheta_0\,x)\,e^{\sin(\vartheta_0
\,x)}\;{\rm
d}x}{\left(\int_{0}^{\tau}\cos^2(\vartheta_0\,x)\,e^{\sin(\vartheta_0
\,x)}\;{\rm d}x\right)^{2}},
\end{equation*}
for any $u>0$.

\bigskip

$\mathbf{Representation\; of\; the\; power}$. We consider the
explicit representation of the power of $\phi_{n}^{*}$ up to
$O(\varepsilon_n^3)$. As we saw
\begin{align*}
\beta_{n}\left(u, \phi_{n}^{*}\right)&= {\cal
N}(a_n)+\frac{\gamma_{3,n}}{3!}\,H_2(a_n)\,n(a_n)
-\frac{\gamma_{4,n}}{4!}\,H_3(a_n)\,n(a_n)-\\&-\frac{\gamma^2_{3,n}}{72}\,H_5(a_n)\,n(a_n)
+O(\varepsilon_{n}^{3}).
\end{align*}
Since $a_n-\Delta=O(\varepsilon_n)$, it can be written as
\begin{align*}
\beta_{n}\left(u, \phi_{n}^{*}\right)&= {\cal
N}(a_n)+\frac{\gamma_{3,n}}{3!}\,H_2(a_n)\,n(a_n)
-\frac{\gamma_{4,n}}{4!}\,H_3(\Delta)\,n(\Delta)-\\&-\frac{\gamma^2_{3,n}}{72}\,H_5(\Delta)\,n(\Delta)
+O(\varepsilon_{n}^{3}).
\end{align*}
Remind that $\Delta=u-z_\alpha$. Indeed we have (see Appendix)
\begin{align*}
a_n-\Delta&=\frac{u^2}{2}\,\varphi_n^3\,I(1,1,1)
+\frac{u^3}{6}\,\varphi_n^4\,I(1,1,0,1)+\\&+\frac{1-z_{\alpha}^2-3u\,\Delta}{6}
\,\varphi_n^3\,I(2,3,0)+\\&+\frac{9u^2\,\Delta-2u(1-z_\alpha^2)}{24}
\,\varphi_n^6\,I(2,3,0)^2-\frac{u^2\,\Delta}{4}\,\varphi_n^4\,I(2,2,1)-\\&-
\frac{u^3}{4}\,\varphi_n^6\,I(2,3,0)\,I(1,1,1)+O(\varepsilon_n^3).
\end{align*}
Taylor expansion yields
$$
{\cal N}(a_n)={\cal
N}(\Delta)+(a_n-\Delta)n(\Delta)-\frac{(a_n-\Delta)^2\Delta\,
n(\Delta)}{2}+O(\varepsilon_n^3)
$$
and
\begin{align*}
\gamma_{3,n}\,H_2(a_n)\,n(a_n)&=\gamma_{3,n}\,H_2(\Delta)\,n(\Delta)+\\&+
\Delta\left(2-H_2(\Delta)\right)\,n(\Delta)\,\gamma_{3,n}\left(a_n-\Delta\right)+O(\varepsilon_n^3).
\end{align*}
Note that
$$
(a_n-\Delta)^2=\left(\frac{u^2}{2}\,\varphi_n^3\,I(1,1,1)
+\frac{1-z_{\alpha}^2-3u\,\Delta}{6}
\,\varphi_n^3\,I(2,3,0)\right)^{2}+O(\varepsilon_n^3)
$$
and
$$
\gamma_{3,n}\left(a_n-\Delta\right)=\frac{u^2}{2}\,\varphi_n^6\,I(2,3,0)\,I(1,1,1)
+\frac{1-z_{\alpha}^2-3u\,\Delta}{6}
\,\varphi_n^6\,I(2,3,0)^{2}+O(\varepsilon_n^3).
$$
Therefore we arrive at the following form
$$
\beta_{n}\left(u, \phi_{n}^{*}\right)= {\cal
N}(\Delta)+n(\Delta)\,\left(r_1(n)+r_2(n)\right)+O(\varepsilon_n^3),
$$
where the terms $r_i(n)=O(\varepsilon_n^i),\; i=1,2,$
\begin{align*}
&r_1(n)=\frac{u^2}{2}\,\varphi_n^3\,I(1,1,1)
+\frac{u\,z_\alpha-2u^2}{6} \,\gamma_{3,n}\\
&r_2(n)=\frac{u^3}{6}\,\varphi_n^4\,I(1,1,0,1)+\frac{9u^2\,\Delta-2u(1-z_\alpha^2)}{24}
\,\gamma_{3,n}^2-\\&-\frac{u^2\,\Delta}{4}\,\varphi_n^4\,I(2,2,1)-
\frac{u^3}{4}\,\varphi_n^3\,\gamma_{3,n}\,I(1,1,1)-\\&-\frac{\Delta}{2}\left(\frac{u^2}{2}\,\varphi_n^3\,I(1,1,1)
+\frac{1-z_{\alpha}^2-3u\,\Delta}{6} \,\gamma_{3,n}\right)^{2}+\\&+
\frac{\Delta\left(2-H_2(\Delta)\right)}{6}\left(\frac{u^2}{2}\,\varphi_n^3\,\gamma_{3,n}\,I(1,1,1)
+\frac{1-z_{\alpha}^2-3u\,\Delta}{6}\, \gamma_{3,n}^{2}\right)-\\&
-\frac{\gamma_{4,n}}{4!}\,H_3(\Delta)\,n(\Delta)-\frac{\gamma^2_{3,n}}{72}\,H_5(\Delta)\,n(\Delta).
\end{align*}
For the notations $I$ above,  see \eqref{not} in the Appendix.
Especially
$$
\varphi_n^3\,I(2,3,0)=\varphi_n^3\int_{\mathbb{A}_{n}}
\frac{\dot{S}\left(\vartheta_{0},
x\right)^{3}}{S\left(\vartheta_{0}, x \right)^{2}}{\rm
d}x=\gamma_{3,n}.
$$
Note also that $r_1(n)=Q_n(u)$ given in the equation
\eqref{power21}.

\section{Appendix}
In order to obtain \eqref{A-a} and \eqref{gamma} we use the
expansion of $S(\vartheta_u,x)$ about $\vartheta_0$ and the
following Taylor expansions which precise the terms up to order
$O(\varepsilon_n^3)$. Letting
$$
t=\frac{S(\vartheta_u,x)-S(\vartheta_0,x)}{S(\vartheta_0,x)},\quad
$$
we have
\begin{align*}
\ln(1+t)&=t-{\frac {1}{2}}{t}^{2}+{\frac {1}{3}}{t}^{3}-{\frac
{1}{4}}{t}^{4}+O
 \left( {t}^{5} \right)\\
 \ln^2(1+t)&={t}^{2}-{t}^{3}+{\frac {11}{12}}{t}^{4}+O \left( {t}^{5}
 \right)\\
\ln^3(1+t)&={t}^{3}-{\frac {3}{2}}{t}^{4}+O \left( {t}^{5} \right)\\
\ln^4(1+t)&={t}^{4}+O \left( {t}^{5} \right).
\end{align*}
Note that the conditions $\mathcal{D}_3$ imply that for the terms of
order $O \left( {t}^{5} \right)$ we have
$$
\int_{\mathbb{A}_n}O \left( {t}^{5} \right)\;{\rm
d}x=O(\varepsilon_n^3).
$$
We use also the binomial expansions
\begin{align*}
(1+s)^{\frac{1}{2}}&=1+{\frac {1}{2}}s-{\frac {1}{8}}{s}^{2}+O
\left(
{s}^{3} \right)\\
(1+s)^{-\frac{1}{2}}&=1-{\frac {1}{2}}s+{\frac {3}{8}}{s}^{2}+O
\left( {s}^{3} \right)\\
(1+s)^{-\frac{3}{2}}&=1-{\frac {3}{2}}s+{\frac {15}{8}}{s}^{2}+O
\left( {s}^{3} \right),
\end{align*}
for different values of $s=s_n(u)=O(\varepsilon_n)$ related to
$\sigma ^2_n(u), \sigma^2_n$ and $\eta^2_n$. We introduce the
following notations:
\begin{equation}
\label{not} I(r_0,r_1,r_2,r_3)=\int_{\mathbb{A}_{n}}
\frac{\dot{S}\left(\vartheta_{0},
x\right)^{r_1}\ddot{S}\left(\vartheta_{0},
x\right)^{r_2}S^{(3)}\left(\vartheta_{0},
x\right)^{r_2}}{S\left(\vartheta_{0}, x \right)^{r_0}}{\rm d}x
\end{equation}
for nonnegative integers $r_0, r_1, r_2$ and
$r_3\in\left\{0,1\right\}$. When $r_3=0$ we write $I(r_0,r_1,r_2)$
instead of $I(r_0,r_1,r_2,0)$.

\noindent Remind that
$$
A_n=\frac{\mu_n(u)-b_n(u)}{\sigma_n(u)} \qquad
a_n=\frac{m_n(u)-c_n}{\eta_n}.
$$
\noindent These give:
\begin{align*}
\mu_n(u)&=\frac{u^2}{2}+\frac{u^4}{8}\,\varphi_n^4\,I(1,0,2)+\frac{u^3}{2}\,\varphi_n^3\,I(1,1,1)
+\frac{u^4}{6}\,\varphi_n^4\,I(1,1,0,1)-\\&-\frac{u^3}{6}\,\varphi_n^3\,I(2,3,0)
-\frac{u^4}{4}\,\varphi_n^4\,I(2,2,1)+\frac{u^4}{12}\,\varphi_n^4\,I(3,4,0)+O(\varepsilon_n^3)\\
\sigma_n(u)^{-1}&=u^{-1}-\frac{u}{8}\,\varphi_n^4\,I(1,0,2)-\frac{1}{2}\,\varphi_n^3\,I(1,1,1)-
\frac{u}{6}\,\varphi_n^4\,I(1,1,0,1)+\\&+\frac{u}{24}\,\varphi_n^4\,I(3,4,0)+
\frac{3u}{8}\,\varphi_n^6\,I(1,1,1)^{2}+O(\varepsilon_n^3)\\
b_n(u)&=-\frac{u^2}{2}+u\,z_\alpha-\frac{u^4-u^3\,z_{\alpha}}{8}\,\varphi_n^4\,I(1,0,2)
-\frac{u^3-u^2\,z_{\alpha}}{2}\,\varphi_n^3\,I(1,1,1)-\\&-\frac{u^4-u^3\,z_{\alpha}}{6}\,\varphi_n^4\,I(1,1,0,1)+
\frac{2u^3-u(1-z_\alpha^2)-3u^2\,z_{\alpha}}{6}\,\varphi_n^3\,I(2,3,0)+\\&+
\frac{2u^4-u^2(1-z_\alpha^2)-3u^3\,z_{\alpha}}{4}\,\varphi_n^4\,I(2,2,1)-\\&-
\frac{6u^4-6u^2(1-z_\alpha^2)-11u^3\,z_{\alpha}}{24}\,\varphi_n^4\,I(3,4,0)-\\&-
\frac{4u^2(1-z_\alpha^2)+3u^3\,z_{\alpha}}{24}\,\varphi_n^6\,I(2,3,0)^2+\\&+
\frac{2u^2(1-z_\alpha^2)+3u^3\,z_{\alpha}}{12}\,\varphi_n^6\,I(2,3,0)\,I(1,1,1)
-\frac{u^3\,z_\alpha}{8}\,\varphi_n^6\,I(1,1,1)^2+O(\varepsilon_n^3).
\end{align*}
Hence we arrive at
\begin{align*}
A_n&=u-z_\alpha+\frac{u^3}{8}\,\varphi_n^4\,I(1,0,2)+\frac{u^2}{2}\,\varphi_n^3\,I(1,1,1)
+\frac{u^3}{6}\,\varphi_n^4\,I(1,1,0,1)-\\&-
\frac{3u^2-(1-z_\alpha^2)-3u\,z_{\alpha}}{6}\,\varphi_n^3\,I(2,3,0)-\\&-
\frac{3u^3-6u(1-z_\alpha^2)-12u^2\,z_{\alpha}}{24}\,\varphi_n^4\,I(3,4,0)+\\&+
\frac{4u(1-z_\alpha^2)+3u^2\,z_{\alpha}}{24}\,\varphi_n^6\,I(2,3,0)^2+\\&+
\frac{3u^3-3u(1-z_\alpha^2)-6u^2\,z_{\alpha}}{12}\,\varphi_n^6\,I(2,3,0)\,I(1,1,1)
-\\&-\frac{u^3\,z_\alpha}{8}\,\varphi_n^6\,I(1,1,1)^2+O(\varepsilon_n^3).
\end{align*}
For the term $a_n$ we have
\begin{align*}
m_n(u)&=u+\frac{u^2}{2}\,\varphi_n^3\,I(1,1,1)
+\frac{u^3}{6}\,\varphi_n^4\,I(1,1,0,1)+O(\varepsilon_n^3),\\
\eta^{-1}_n&=1-\frac{u}{2}\,\varphi_n^3\,I(2,3,0)-\frac{u^2}{4}\,\varphi_n^4\,I(2,2,1)+\frac{3u^2}{8}
\,\varphi_n^6\,I(2,3,0)^2+O(\varepsilon_n^3).
\end{align*}
Therefore we can write
\begin{align*}
a_n&=u-z_{\alpha}+\frac{u^2}{2}\,\varphi_n^3\,I(1,1,1)
+\frac{u^3}{6}\,\varphi_n^4\,I(1,1,0,1)+\\&+\frac{1-z_{\alpha}^2-3u(u-z_\alpha)}{6}
\,\varphi_n^3\,I(2,3,0)+\\&+\frac{9u^2(u-z_{\alpha})-2u(1-z_\alpha^2)}{24}
\,\varphi_n^6\,I(2,3,0)^2-\frac{u^2(u-z_\alpha)}{4}\,\varphi_n^4\,I(2,2,1)-\\&-
\frac{u^3}{4}\,\varphi_n^6\,I(2,3,0)\,I(1,1,1)+O(\varepsilon_n^3).
\end{align*}
Hence one can obtain $A_n-a_n$ as in \eqref{A-a}. Using the
necessary Taylor expansions we get
\begin{align*}
\gamma_{3,n}(u)&=\varphi_n^3\,I(2,3,0)+u\,\varphi_n^4\,I(3,4,0)-\frac{3u}{2}
\,\varphi_n^6\,I(2,3,0)^2+O(\varepsilon_n^3),\\
\gamma'_{3,n}(u)&=\varphi_n^3\,I(2,3,0)-\frac{u}{2}\,\varphi_n^4\,I(3,4,0)+\frac{3u}{2}
\,\varphi_n^4\,I(2,2,1)-\\&-
\frac{3u}{2}\,\varphi_n^6\,I(2,3,0)\,I(1,1,1)+O(\varepsilon_n^3),
\end{align*}
which prove \eqref{gamma}.

 \end{document}